\documentclass{jgcc}
\pdfoutput=1
\usepackage{lastpage}
\jgccdoi{17}{1}{1}{13459} 
\jgccheading{}{\pageref{LastPage}}{}{}{Apr.~24,~2024}{Feb.~20,~2025}{}   

\usepackage{amsfonts, amsmath, amssymb, amscd, amsthm, graphicx, setspace, enumitem, color}
\usepackage{hyperref}
\usepackage{pinlabel}
\usepackage{mathtools}
\usepackage{bm}
\usepackage[normalem]{ulem}
\usepackage{graphicx}

\makeatletter
\newtheorem*{rep@theorem}{\rep@title}
\newcommand{\newreptheorem}[2]{%
\newenvironment{rep#1}[1]{%
 \def\rep@title{#2 \ref{##1}}%
 \begin{rep@theorem}}%
 {\end{rep@theorem}}}
\makeatother

\newreptheorem{theorem}{Theorem}

\newtheoremstyle{exstyle}{}{}{}{}{\bfseries}{ }{5pt}{}
\theoremstyle{exstyle}

\theoremstyle{definition}

\newcommand{\ie}{\textit{i.e.},~}

\newcommand{\st}{\;\colon\;}
\newcommand{\bs}{\backslash}
\newcommand{\defn}[1]{\textcolor{blue}{\emph{#1}}}
\newcommand{\set}[2]{\left\{#1\st#2\right\}} 
\newcommand{\gen}[1]{\left<#1\right>}
\newcommand{\genset}[2]{\gen{#1\st#2}}
\DeclarePairedDelimiter\order{\lvert}{\rvert}%
\newcommand{\iso}{\cong}

\newcommand{\orig}{o}
\DeclareMathOperator{\Mov}{Mov}
\DeclareMathOperator{\Fix}{Fix}
\DeclareMathOperator{\vspan}{span}

\renewcommand{\cc}[1]{C\left(#1\right)}

\newcommand{\R}{{\mathbb R}}
\newcommand{\N}{{\mathbb N}}
\newcommand{\Z}{{\mathbb Z}}

\newcommand{\cA }{\mathcal A}

\DeclareMathOperator{\vrank}{rank}

\begin{document}

\title{Conjugacy Class Growth in Virtually Abelian Groups}

\author[A.~Dermenjian]{Aram Dermenjian}	
\address{The University of Manchester, Manchester, UK, and Heilbronn Institute for Mathematical Research, Bristol, UK} 
\email{aram.dermenjian.math@gmail.com} 

\author[A. Evetts]{Alex Evetts}
\address{The University of Manchester, Manchester, UK, and Heilbronn Institute for Mathematical Research, Bristol, UK}
\email{alex.evetts@manchester.ac.uk}

 






\begin{abstract}
    We study the growth of conjugacy classes in finitely generated virtually abelian groups. That is, the number of elements in the ball of radius $n$ in the Cayley graph which intersect a fixed conjugacy class. In the class of virtually abelian groups, we prove that this function is always asymptotically equivalent to a polynomial. Furthermore, we show that in any affine Coxeter group, the degree of polynomial growth of a conjugacy class is equivalent to the reflection length of any element of that class.
\end{abstract}

\maketitle


\section{Introduction}
Given a finitely generated group $G$ generated by some set $S$, every element in $G$ can be represented by some word over the alphabet $S \cup S^{-1}$.
For an element $g \in G$, we define the length $\ell_S(g)$ of $g$ to be the minimal length of all words which represent $g$ with respect to the generating set $S$.
This naturally leads to the study of the growth function of $G$, which counts the number of elements of length at most $n$. The asymptotic properties of growth functions of groups have been extensively studied with a highlight being Gromov's celebrated theorem \cite{Gromov} characterising groups of polynomial growth as the virtually nilpotent groups.

The growth function can be generalised by restricting our attention to subsets of interest.
In particular, taking some subset $U$ of $G$, we let $\beta_{U, S}(n)$ be the number of elements in $U$ whose length is at most $n$ and look at how this number changes as $n$ increases.
The growth of subgroups of groups has been studied by a number of authors \cite{CantrellSharp, DavisOlshanskii, Schesler} but much less attention has been paid to other interesting subsets.

In this article we focus our attention on conjugacy classes in virtually abelian groups.
In particular, we take a conjugacy class $C$ and we look at the growth function of $C$ as a subset of $G$.
To this end, we prove the following theorem.
\begin{reptheorem}{thm:conj_class_poly}
    Let $G$ be a finitely generated virtually abelian group. Then every conjugacy class has polynomial growth.
\end{reptheorem}
We emphasise that this is not to be confused with the \emph{conjugacy growth function}, which counts the number of conjugacy classes intersecting the ball of radius $n$.
The conjugacy growth function in virtually abelian groups was studied in~\cite{Evetts}. Combinatorial properties including growth of other subsets of virtually abelian groups are explored in \cite{RationalVAb}, and as noted in Remark \ref{rem:altproof}, Theorem \ref{thm:conj_class_poly} also follows from general results in that article. However, the method of proof we use here is much more direct, and as we will see, allows for the degree of polynomial growth to be explicitly calculated.

To the authors' knowledge, the first appearance of the growth of a conjugacy class is in \cite{ParkPaul}, where it appears as a special case of the `orbital counting problem in conjugacy classes'. The authors of that paper are interested in examples of groups where all non-trivial conjugacy classes have equivalent growth functions, as well as the \emph{exponential growth rate} of the function. They show that in the cases of hyperbolic groups, and the integer Heisenberg group, all infinite conjugacy classes have the same growth. In the present paper we show that this is not the case in general for virtually abelian groups (see Proposition~\ref{prop:allgrowths}), while providing an example of where it does hold (Example~\ref{ex:kleinbottle}).

As an example of the growth function of conjugacy classes in virtually abelian groups, we consider the affine Coxeter groups.
It turns out that the polynomial growth is directly related to the reflection length.
For an element $w$ in an affine Coxeter group $W$ with simple reflections $S$, we decompose $w$ into a translation part $t$ and a finite part $u$ where $w = tu$.
The reflection length $\ell_R(w)$ of $w$ is the length with regards to the set of \emph{all} reflections instead of just the simple reflections.
This leads us to our second main theorem of this article.
\begin{reptheorem}{thm:poly_growth_Coxeter}
    Let $W$ be an affine Coxeter group and $w = tu$ an arbitrary element in $W$ decomposed as a translation part and a finite part.
    Then the growth of the conjugacy class with respect to the generating set $S$ is equivalent to $n^{\ell_R(u)}$.
\end{reptheorem}

This article is organised in the following way.
We start by giving preliminary definitions and basic results in group theory and growth in Section~\ref{sec:definitions}.
In Section~\ref{sec:coxeter_groups}, we give background definitions and results on affine Coxeter groups.
The affine Coxeter groups are a particularly nice class of virtually abelian groups as they have a rich geometric structure to them.
Section~\ref{sec:virtually_abelian_groups} is then devoted to showing our two main results: (1) that the growth function of conjugacy classes in virtually abelian groups is polynomial and (2) that the growth function of the conjugacy class of an element in an affine Coxeter group is polynomial with degree equal to the reflection length of its finite part.

\section{Definitions and basic results}
\label{sec:definitions}
We begin by recalling some standard definitions in group theory and growth functions. The interested reader is referred to \cite{Loh} and \cite{Mann} for a more thorough introduction to the growth of groups.

\subsection{Growth in groups}
\label{ssec:growth_in_groups}
We start by giving some basic notions on words representing elements of a group and the growth rate of a set in relation to a group.
\begin{defi}\label{def:growth}
    Let $G$ be a finitely generated group and $S$ a choice of finite generating set. We let $\left\{ S \cup S^{-1} \right\}^{\star}$ denote the words over the alphabet $S$ and its inverses, including the empty word.
    Given a word $w \in \left\{ S \cup S^{-1} \right\}^{\star}$ we let $\order{w}$ denote the length of the word $w$ and we write $w =_G g$ if $w$ represents $g$ in the group $G$.
    Many words can represent the same element, but there will always be some word of minimal length representing a given $g$.
    The \defn{word length} of an element $g\in G$ with respect to $S$ is then defined to be this minimal length
    \begin{equation*}
        \ell_S(g)=\min\set{\order{w}}{w\in\{S\cup S^{-1}\}^*,w=_Gg}.
    \end{equation*}
    If there is no ambiguity in the generating set, we will write $\ell(g)$ as a shorthand for $\ell_S(g)$.
    The length of elements can be extended to the \defn{word metric} on $G$ via the distance function given by $d_{G,S}(g,h)=\ell_S(g^{-1}h)$.
    We write $B_{G,S}(n)=\set{g\in G}{ \ell_S(g)\leq n}$ for the metric $n$-ball in $G$ with respect to $S$.
    If there is no ambiguity in the group or generating set, we will use $B(n)$ to represent $B_{G, S}(n)$.
    Let \(U\) be any subset of \(G\). We define the \defn{growth of $U$ in $G$} as follows:
    \begin{equation*}
        \beta_{U,S}(n) \coloneqq \order{\set{g\in U}{\ell_{S}(g)\leq n}} = \order{U\cap B_{G,S}(n)}.
    \end{equation*}
    As before, if there is no ambiguity in the set $S$, the group $G$ or the subset $U$, instead of $\beta_{U, S}(n)$ we will write $\beta_U(n)$. When $U=G$ this is the familiar notion of the (cumulative) growth function of a group.
\end{defi}

The following lemma shows that the number of elements in the $n$-ball is the same, regardless of the choice of centre point.
\begin{lem}\label{lem:nballs}
    Fix any element $h\in G$. Then $\order{\set{g\in G}{d_{G,S}(h,g)\leq n}}=\beta_{G}(n)$.
\end{lem}
\begin{proof}
    For $h,x,y\in G$ we have $d_{G,S}(x,y)=\ell_S(x^{-1}y)=\ell_S(x^{-1}h^{-1}hy)=d_{G,S}(hx,hy)$, and so the word metric is invariant under left multiplication. In other words $d_{G,S}(1,x)\leq n$ if and only if $d_{G,S}(h,hx)\leq n$ and we have a length-preserving bijection between the elements of the $n$-ball centred at $1$ and the elements of the $n$-ball centred at $h$.
\end{proof}

Note that the growth of a subgroup $U$ of $G$ is often referred to as \emph{relative growth}, since the subgroup may itself be finitely generated and therefore have an intrinsic notion of growth. In this paper, we study subsets that are not necessarily subgroups and so we drop the term `relative' since the natural growth is the one inherited from the ambient group.

We recall the usual notion of equivalence of growth functions. Let $f,g\colon\N\rightarrow\N$ be two functions. We write $f\preccurlyeq g$ if there exists $\lambda\geq1$ such that
    \[f(n)\leq\lambda g(\lambda n+\lambda)+\lambda\] for all $n\in\N$. If $f\preccurlyeq g$ and $g\preccurlyeq f$ then we write $f\sim g$ and say that the functions are \defn{equivalent}. Note that this defines an equivalence relation. We write $f\prec g$ if $f\preccurlyeq g$ but it is not the case that \(f\sim g\).

It is a standard fact that, up to equivalence, the growth of a group does not depend on the choice of finite generating set (since every element written as a word of one generating set can be written as a word of bounded length in another generating set). The following lemma is easily proved in the same way.
\begin{lem}\label{lem:relgrowthwelldefined}
    Up to equivalence, the growth of a subset $U$ in $G$ does not depend on the choice of generating set for $G$.
\end{lem}

We say that a subset $U$ has \defn{polynomial growth} in $G$, if $\beta_{U}(n)\sim n^d$ for some positive integer $d$.
We will also make use of the following elementary fact.
\begin{lem}\label{lem:union}
    Let \(U\) and \(V\) be subsets of a group \(G\). Then the growth of \(U\cup V\) in $G$ is equivalent to the maximum (with respect to $\preccurlyeq$) of the growths of \(U\) and \(V\).
\end{lem}
\begin{proof}
    Without loss of generality, suppose that $\beta_U\preccurlyeq\beta_V$. Then we have $\beta_V\preccurlyeq\beta_{U\cup V}\preccurlyeq\beta_{U}+\beta_V\sim\beta_V$,
    and so $\beta_V\sim\beta_{U\cup V}$.
\end{proof}

In Lemma~\ref{lem:relgrowthwelldefined} we saw that changing the generating set does not change the asymptotics of the growth function. Furthermore, just as with the standard growth function, the growth of a subset is robust enough to be an invariant of the usual large-scale geometry of a group.

Recall the following definition of a quasi-isometric embedding of metric spaces.
    Let $(X,d)$ and $(Y,d')$ be metric spaces. A function $f\colon X\to Y$ is a \defn{quasi-isometric embedding} if there exists a constant $\lambda\geq1$ such that for every $x,y\in X$ we have
    \[\frac{1}{\lambda}d(x,y)-\lambda \leq d'(f(x),f(y)) \leq \lambda d(x,y)+\lambda.\]
    
    Any subset of a finitely generated group can be given the structure of a metric space via Definition~\ref{def:growth} and therefore we may consider quasi-isometric embeddings between subsets of finitely generated groups. It is a standard result (see, for example, \cite[Proposition 6.2.4]{Loh}) that quasi-isometric groups have equivalent growth functions. The next result can be proved in the same way, using Lemma~\ref{lem:nballs}.

\begin{prop}\label{prop:qi}
    Let $f\colon G\to H$ be a quasi-isometric embedding of finitely generated groups. Let $U\subset G$. Then the growth of $U$ in $G$ is equivalent to the growth of $f(U)$ in $H$, i.e. $\beta_{U}\sim\beta_{f(U)}$.
\end{prop}
A special case of Proposition \ref{prop:qi} is the following.
\begin{cor}\label{cor:translatedgrowth}
    Let $U$ be any subset of a group $G$. Then for any $g\in G$, the subset $Ug=\set{ug}{u\in U}$ has equivalent growth to $U$.
\end{cor}

The subject of this article is \defn{conjugacy class growth}, the growth function $\beta_C(n)$ for a fixed conjugacy class $C$ of some group. By Lemma \ref{lem:relgrowthwelldefined}, this does not depend on the choise of generating set.




The main class of groups in which we will study conjugacy class growth are virtually abelian groups.
A well-known example of such groups are affine Coxeter groups which we define in Section~\ref{sec:coxeter_groups}.

\section{Virtually abelian groups}
\label{sec:virtually_abelian_groups}
Recall that, given some property $P$, a group is said to be \defn{virtually $P$} if it has a finite-index subgroup with property $P$. By a famous result of Bass and Guivarc'h~\cite{Bass}, all finitely generated virtually nilpotent groups have polynomial growth. Therefore the conjugacy classes of such a group have growth bounded above by a polynomial. In the case of virtually abelian groups, we show that the conjugacy class growth function is always \emph{equivalent} to a polynomial, thus excluding the possibility of growth functions that lie between polynomials (such as \(n^{3/2}\), \(n\log n\)).

\subsection{Virtually abelian groups}
We start by looking at the conjugacy class growth of virtually abelian groups in general.
\begin{rem}\label{rem:normalcore}
    If $G$ is virtually abelian, it is a standard fact that $G$ contains a finite-index \emph{normal} abelian subgroup (for example, the normal core of any finite-index abelian subgroup, the intersection of all conjugates of that subgroup, see for example \cite[Chapter 2]{Mann}). 
\end{rem}
The proof of the following proposition is a standard exercise.
    \begin{prop}\label{prop:abeliangrowth}
        The growth function of the free abelian group $\Z^r$ is equivalent to $n^r$, for all positive integers $r$.
    \end{prop}

\begin{defi}
    Let $G$ be a group with a subgroup $H$. For any $g\in G$, define
    \[[H,g]=\langle [h,g]\mid h\in G\rangle.\]
\end{defi}
\begin{lem}\label{lem:commlength}
    Let $G$ be a group with a normal abelian subgroup $H$. Then for any $g\in G$,
    \[
        [H,g] = \set{[h,g]}{h\in H}.
    \]
    In other words, every element of $[H,g]$ has commutator length equal to one.
\end{lem}
\begin{proof}
    Let $h,h'\in H$. By the normality of $H$ we have $gh^{-1}g^{-1}\in H$, and so $h'$ and $gh^{-1}g^{-1}$ commute, and thus we have 
    \begin{align*}
        [h,g][h',g] = hgh^{-1}g^{-1} h'gh'^{-1}g^{-1} = hh'gh^{-1}g^{-1}gh'^{-1}g^{-1}.
    \end{align*}
    Cancelling the $g^{-1}g$ and noting that $h^{-1}h'^{-1}=(h'h)^{-1}=(hh')^{-1}$ as $H$ is abelian gives 
    \begin{align*}
        [h,g][h',g]= hh'gh^{-1}h'^{-1}g^{-1} = [hh',g],
    \end{align*}
    which finishes the proof.
\end{proof}

\begin{rem}\label{rem:MST1}
    The following characterisation of conjugacy classes is similar to that used in \cite[Section 5]{Evetts}. Furthermore, after a preprint version of the current paper was published on arXiv, a preprint of Mili\'cevi\'c, Schwer and Thomas \cite{MST} appeared, using a similar description of conjugacy classes, in the case of affine Coxeter groups, with so-called `Mod-sets' in place of subgroups of commutators (see Theorem 2.12 in that paper). See also Remark \ref{rem:MST2}.
\end{rem}
For any group element $g$, write $\cc{g}=\set{hgh^{-1}}{h\in G}$ for its conjugacy class.
\begin{lem}\label{lem:vabconj}
    Let $G$ be a virtually abelian group, and let $H$ be a finite-index normal abelian subgroup of $G$. Let $U$ be any set of coset representatives for $H$. Then the conjugacy class of any $g\in G$ is a union of cosets of the form
    \begin{align*}
        \cc{g} &= \bigcup_{v\in U} [H,vuv^{-1}]vgv^{-1},\\
    \end{align*}
    for a fixed $u\in U$.
\end{lem}
\begin{proof}
    We have \(g=hu\) for some \(h\in H\), \(u\in U\).
    The conjugacy class of \(g\) is then
    \begin{align*}
        \cc{g}=\cc{hu} &=\set{xv hu v^{-1}x^{-1}}{x\in H,\,v\in U}\\
        &= \set{xvhv^{-1}vuv^{-1}x^{-1}}{x\in H,\,v\in U}\\
        &= \bigcup_{v\in U}\set{xvhv^{-1}vuv^{-1}x^{-1}vu^{-1}v^{-1}}{ x\in H}vuv^{-1}\\
        &= \bigcup_{v\in U}\set{xvuv^{-1}x^{-1}vu^{-1}v^{-1}}{ x\in H}vhv^{-1}vuv^{-1}\\
        &= \bigcup_{v\in U}\set{[x,vuv^{-1}]}{ x\in H}vhuv^{-1}\\
        &=\bigcup_{v\in U}[H,vuv^{-1}]vgv^{-1},
    \end{align*}
    where we have used normality of $H$ in $G$ in the fifth equality and Lemma~\ref{lem:commlength} in the last equality.
\end{proof}

{Before we prove the main theorem, we give a short lemma.}
\begin{lem}
    \label{lem:bilipschitz_emb}
    Let $G$ be a finitely generated free abelian group and let $H$ be any subgroup of $G$.
    The inclusion map $i \colon H \to G$ is a quasi-isometric embedding.
\end{lem}
\begin{proof}
    Any subgroup of a free abelian group is itself free abelian, of equal or smaller rank, so let $A=\{a_1,\ldots,a_k\}$ be a basis for $H$. Then there exists a basis $B=\{b_1,\ldots,b_k,\ldots,b_l\}$ for $G$ such that for each $1\leq j\leq k$ there exists a positive integer $\lambda_j$ with $i(a_j)=\lambda_jb_j$. In particular, write $i(a_j)=(x_1,\ldots,x_d)\in\Z^d$ in terms of the standard basis for $G$. Then $\lambda_j=\gcd(\order{x_1},\ldots,\order{x_d})$, and $b_j=(x_1/\lambda_j,\ldots,x_d/\lambda_j)$.
    
    Let $h\in H$. Then, writing the group operation additively, there is a unique way to write $h=\sum_{j=1}^k\mu_ja_j$, with $\mu_j\in\Z$, and so $\ell_A(h)=\sum_{j=1}^k\order{\mu_j}$.
    Applying the inclusion map gives $i(h)=\sum_{j=1}^k\mu_j\lambda_ja_j$, again a sum of basis elements, so $\ell_B(i(h))=\sum_{j=1}^k\order{\mu_j\lambda_j}$. Combining these equalities, we have
    \[
        \ell_A(h)\leq \ell_B(i(h))\leq\lambda \ell_A(h),
    \]
    where $\lambda=\max_j\lambda_j$, which implies that the map $i$ is a quasi-isometric embedding.
\end{proof}

\begin{thm}
    \label{thm:conj_class_poly}
    Let \(G\) be a finitely generated virtually abelian group. Then every conjugacy class has polynomial growth.
\end{thm}
\begin{proof}
    We may assume that $G$ has a finite-index free abelian normal subgroup $H$ (by first passing to a free abelian subgroup, and then to its normal core as per Remark~\ref{rem:normalcore}).
    Then by Lemmas~\ref{lem:vabconj} and~\ref{lem:union} and Corollary~\ref{cor:translatedgrowth}, it is enough to show that the growth of any subgroup $[H,g]$ in $G$ is polynomial.
    
    By normality of $H$, $[H,g]$ is contained in $H$ and is therefore free abelian of finite rank, and so by Proposition \ref{prop:abeliangrowth} it has polynomial growth (of degree equal to its rank).
    By Lemma~\ref{lem:bilipschitz_emb}, we have that $[H,g]$ is quasi-isometrically embedded into $H$.
    Since $H$ is a finite-index subgroup of $G$, it is quasi-isometrically embedded into $G$, and so $[H,g]$ is also quasi-isometrically embedded into $G$.
    By Proposition~\ref{prop:qi}, the growth of $[H,g]$ in $G$ is equivalent to the growth of $[H,g]$ itself and is therefore polynomial.
\end{proof}

\begin{rem}\label{rem:altproof}
    An alternative and less explicit method of proof would be to observe that Lemma~\ref{lem:vabconj} implies that a conjugacy class is an example of a \defn{rational subset} (the image of a \defn{regular language} over the generating set).
    In \cite{RationalVAb} it is proved that the growth of any rational subset of a virtually abelian group has rational generating function, with respect to any choice of generating set. Any non-decreasing integer-valued sequence with growth in the polynomial range and rational generating function necessarily has growth equivalent to a polynomial (see, for example, \cite{Stanley}).
\end{rem}

The following proposition demonstrates that there exist groups with any given degree of polynomial conjugacy class growth.
\begin{prop}\label{prop:allgrowths}
    For any \(d\in\N\) there exists a (virtually abelian) group \(G\) such that for every natural number $c \in \N$ such that \(c\leq d\), there is a conjugacy class of \(G\) with growth equivalent to \(n^c\).
\end{prop}
\begin{proof}
    Consider the free abelian group \(\Z^d=\langle t_1,\ldots,t_d\rangle\) and extend it by the direct product of $d$ copies of $C_2$, where the conjugation action of the $i$th copy negates the $i$th coordinate of $\Z^d$. In other words, let 
    \begin{align*}
        G &=\Z^d\rtimes(C_2\times\cdots\times C_2)\\
        &=\langle t_1,\ldots,t_d,s_1,\ldots,s_d\mid [t_i,t_j]=[s_i,t_j]=1~\text{for }i\neq j, s_i^2=1, s_it_i=t_i^{-1}s_i \rangle.
   \end{align*}
   So we have 
    \[s_it_js_i=\begin{cases} t_j^{-1} & i=j \\ t_j & i\neq j. \end{cases}\]
The set \(\set{t_1^{k_1}t_2^{k_2}\cdots t_d^{k_d}s_1^{\varepsilon_1}s_2^{\varepsilon_2}\cdots s_d^{\varepsilon_d}}{ k_i\in\Z,\varepsilon_j\in\{0,1\}}\) is a normal form for elements of the group and conjugation by generators behaves as follows:
    \begin{align*}
        s_i t_1^{k_1}t_2^{k_2}\cdots t_d^{k_d}s_1^{\varepsilon_1}s_2^{\varepsilon_2}\cdots s_d^{\varepsilon_d} s_i &= t_1^{k_1}t_2^{k_2}\cdots t_i^{-k_i}\cdots t_d^{k_d}s_1^{\varepsilon_1}s_2^{\varepsilon_2}\cdots s_d^{\varepsilon_d} \\
        t_i t_1^{k_1}t_2^{k_2}\cdots t_d^{k_d}s_j t_i^{-1} &= \begin{cases} t_1^{k_1}t_2^{k_2}\cdots t_d^{k_d}s_j & i\neq j \\ t_1^{k_1}t_2^{k_2}\cdots t_i^{k_i+2}\cdots t_d^{k_d}s_i & i=j\end{cases}.
    \end{align*}
    Therefore for any natural number \(c\leq d\), the conjugacy class of the element \(s_1s_2\cdots s_c\) is the coset \(\langle t_1^2,t_2^2,\ldots,t_c^2\rangle s_1s_2\cdots s_c\) of the free abelian subgroup \(\langle t_1^2,t_2^2,\ldots,t_c^2\rangle\).
    This subgroup has rank \(c\) and therefore it and its cosets (by Corollary \ref{cor:translatedgrowth}) have growth equivalent to \(n^c\). 
\end{proof}

Before we study affine Coxeter groups, we use Lemma~\ref{lem:vabconj} to calculate the conjugacy class growth functions of a virtually abelian group which \emph{does not} split as a semidirect product of the finite-index abelian subgroup.
\begin{exa}\label{ex:kleinbottle}
    Let $K=\langle a,b\mid aba^{-1}=b^{-1}\rangle$, the fundamental group of the Klein bottle, which is a non-split extension of the subgroup \(A=\langle a,b^2\rangle\cong\Z^2\). This subgroup has index $2$ and we choose $\{1,b\}$ as a set of coset representatives.
    If $g\in A$ then both $a$ and $b^2$ commute with $g$.
    Therefore $[g]=\{g,bgb^{-1}\}$ and thus the growth of any conjugacy class contained in $A$ is constant.
    
    If $g\in Ab$ then by Lemma \ref{lem:vabconj} we have \[[g]=[A,b]g \cup [A,bbb^{-1}]bgb^{-1} = [A,b]g \cup [A,b]bgb^{-1}.\] Since $[a,b]=a^2$, we have $[A,b]=\langle a^2\rangle$, which is a cyclic direct factor of $A$, and hence quasi-isometrically embedded and of linear growth. It follows from Proposition \ref{prop:qi}, Corollary \ref{cor:translatedgrowth}, and Lemma \ref{lem:union}, that $\beta_{[g]}(n)$ is linear for conjugacy classes $[g]$ contained in the coset $Ab$.
\end{exa}

\section{Affine Coxeter groups}
\label{sec:coxeter_groups}
In this section we recall the background details needed on affine Coxeter groups.
The reader is assumed to be familiar with Coxeter groups at the level of the book by Humphreys \cite{Humphreys}.
Much of the terminology and notation in this section is taken from \cite{LewisMcCammondPetersenSchwer}.

\subsection{Reflection length}
Given an (arbitrary) Coxeter system $(W, S)$, the set $S$ is called the set of \defn{simple reflections} and $W$ is called the \defn{Coxeter group}.
Recall from Subsection~\ref{ssec:growth_in_groups} that the length function $\ell_S(w)$ gives the minimal length of  words over the alphabet $S$ which represent $w \in W$.
Given a Coxeter system $(W, S)$ we denote by $\Phi$ the associated root system (see \cite{Humphreys} for a definition) with positive roots $\Phi^+$ and negative roots $\Phi^- = -\Phi^+$.

In this paper we make use of another length function which we describe next.
The \defn{set of reflections $R$} is the set of all conjugates to $S$, \ie $R = \bigcup_{w \in W}wSw^{-1}$.
The set of reflections is infinite whenever $W$ is infinite.
As with the simple reflections, any element $w \in W$ can be represented by words over the alphabet $R$.
The \defn{reflection length $\ell_R(w)$} of an element is then the smallest integer $k$ such that $w = r_1 r_2 \ldots r_k$ with $r_i \in R$.
Whenever $w = r_1 \ldots r_k$ such that $\ell_R(w) = k$ then we say that $r_1 \ldots r_k$ is a \defn{reduced expression of $w$ over $R$}.

\begin{exa}
    \label{ex:A2}
    Suppose that $W$ is the type $A_2$ Coxeter group with simple reflections $S = \left\{ s, t \right\}$.
    Recall that being a type $A_2$ Coxeter group implies that we have the following presentation
    \[
        W = \genset{S}{s^2 = t^2 = (st)^3 = 1_W}.
    \]
    The set of all reflections is then given by $R = \left\{ s, t, sts \right\}$.
    The following is a table of the lengths of all elements in $W$:
    \[
        \begin{array}{r || c | c}
            w & \ell(w) & \ell_R(w)\\
            \hline
            1_W & 0 & 0\\
            s & 1 & 1\\
            t & 1 & 1\\
            st & 2 & 2\\
            ts & 2 & 2\\
            sts & 3 & 1
        \end{array}
    \]
    Note that an element might have more than one reduced expression.
    For example, over $R$, $st$ can also be represented as $st = ttst = tsts = stss$.
    Similarly, over $S$, $sts$ can also be represented as $tst$.
\end{exa}

\subsection{Affine Coxeter groups}
\label{ssec:affine_coxeter_groups}
Let $E$ be a Euclidean space with underlying Euclidean vector space $V$.
For any two points $x, y \in E$ recall that there is a unique vector $\lambda \in V$ such that $x + \lambda = y$.
Let $\orig$ denote some (fixed) point in $E$ which we call the \defn{origin}.
Then $E$ and $V$ can be identified by sending each vector $\lambda \in V$ to the point $\orig + \lambda \in E$.

Suppose that $V$ contains a (finite) crystallographic root system $\Phi$ with simple system $\Delta$.
We generate an affine Coxeter group $W$ using $\Phi$ in the following way.
For each $\alpha \in \Phi^+$ and $j \in \Z$, let $H_{\alpha, j} = \set{\lambda \in V}{\gen{\lambda, \alpha} = j}$ be the affine hyperplane in $E$ where $\gen{\cdot,\cdot}$ is the standard inner product.
We let $\cA_W$ denote the set of hyperplanes $H_{\alpha, j}$ for all $\alpha \in \Phi^+$ and $j \in \Z$.
We call $\cA_W$ the \defn{hyperplane arrangement associated to $W$}.
Each affine hyperplane $H_{\alpha, j}$ has a reflection $r_{\alpha, j}$ associated to it which fixes $H_{\alpha, j}$ pointwise.
An \defn{affine Coxeter group $W$} is the group generated by the set of reflections $R = \set{r_{\alpha, j}}{\alpha \in \Phi,\, j \in \Z}$.

To each affine Coxeter group $W$ is associated a finite part $W_0$ which is the finite Coxeter group generated by the restriction of the set $R$ to the set $R_0 = \set{r_{\alpha,0}}{\alpha \in \Phi^+}$.
There is a well-defined projection map $\pi : W \to W_0$ which sends each $r_{\alpha, j}$ to $r_{\alpha, 0}$ which extends to a group homomorphism.
The kernel of this homomorphism is the normal abelian subgroup denoted by $T$, isomorphic to $\Z^d$, whose elements are called \defn{translations} and where $W_0 \iso W / T$ and $d$ is equal to the rank of $W_0$ (and of $W$).
The \defn{simple reflections $S$} of $W$ are the reflections associated to the bounded alcove assocaited to $1_W$, $S = \set{r_{\alpha, 0}}{\alpha \in \Delta} \cup \set{r_{\alpha, 1}}{\alpha \in \left\{\Phi \bs \Delta\right\}}$.

\begin{exa}
    \label{ex:affineA2}
    Let $W$ be the type $\tilde{A}_2$ affine Coxeter group whose underlying finite part is the type $A_2$ Coxeter group from Example~\ref{ex:A2}, \ie $W_0 \iso A_2$.
    The group $W$ is then generated by the simple reflections given by $S = \left\{ r_{\alpha_r, 1}, r_{\alpha_s, 0}, r_{\alpha_t, 0} \right\}$.
    Letting $r = r_{\alpha_r, 1},\, s = r_{\alpha_s, 0},\, t = r_{\alpha_t, 0}$ gives us the presentation
    \[
        W = \genset{S}{r^2 = s^2 = t^2 = (rs)^3 = (st)^3 = (rt)^3 = 1_W}.
    \]
    Then the hyperplane arrangement $\cA_W$ is the set of hyperplanes as given in the following figure.
    \begin{center}
        \includegraphics[scale=0.2]{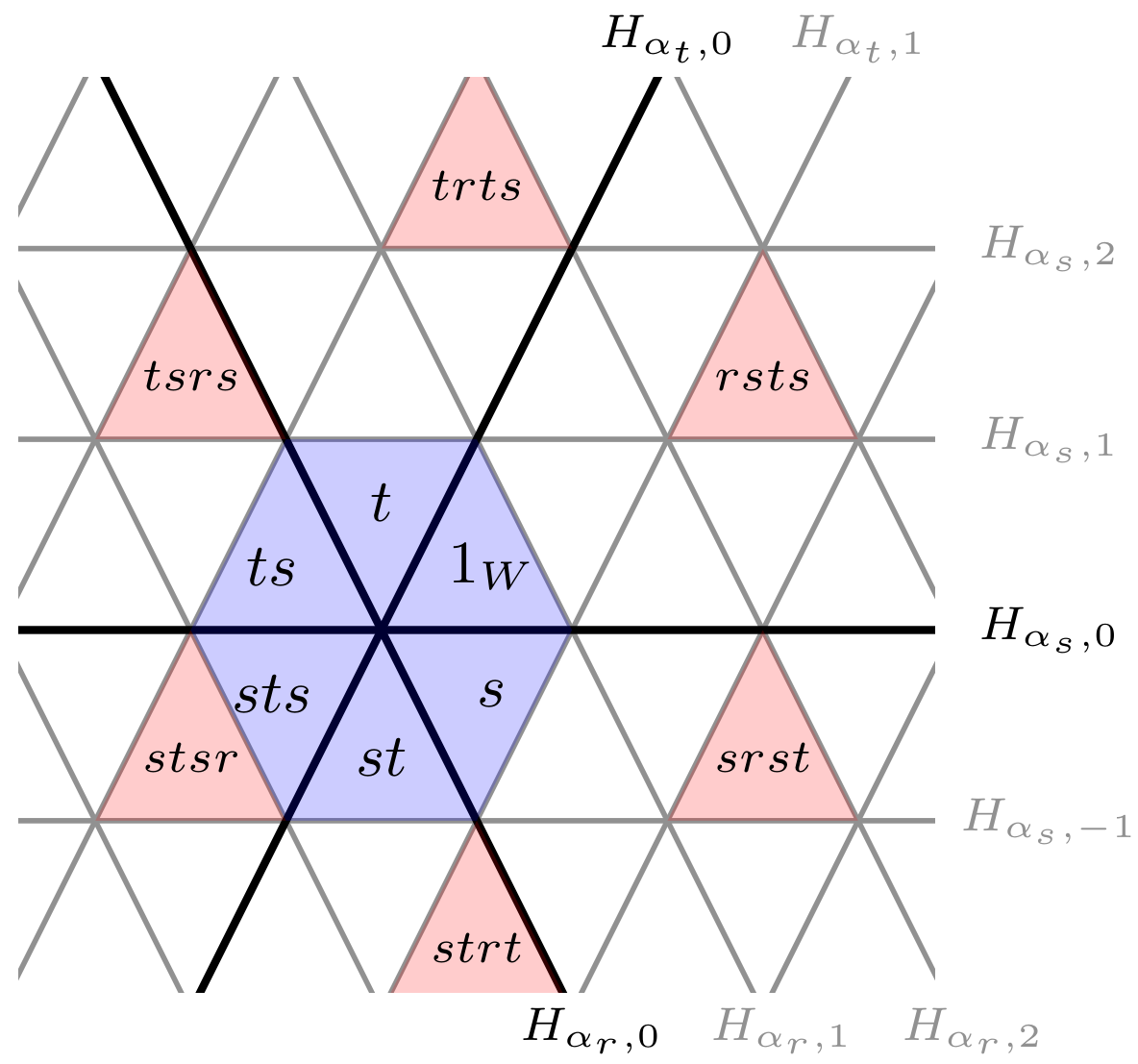}
    \end{center}
    The regions in the figure are labelled by elements from $W$.
    The finite part $W_0$ is shaded in blue and the translations are shaded in red (forming a free abelian subgroup of rank $2$).
\end{exa}

\subsection{Movement}
\label{ssec:movement}
Let $x \in E$ be some point and let $w$ be some Euclidean isometry.
The \defn{motion} of $x$ under $w$ is the vector $\lambda \in V$ such that $w(x) = x + \lambda$.
The \defn{move-set $\Mov(w)$} of $w$ is the collection of all motions of the points in $E$:
\[
    \Mov(w) = \set{\lambda}{w(x) = x + \lambda\text{ for some } x \in E} \subseteq V.
\]
It turns out that $\Mov(w)$ is an affine subspace of $V$ (see \cite[Proposition 3.2]{BradyMcCammond}).
The \defn{fixed space $\Fix(w)$} of $w$ is the set of points $x \in E$ such that $w(x) = x$.
It is an (affine) subspace of $E$ whenever it is nonempty.

An element $w \in W$ is called \defn{elliptic} if its fixed space is nonempty.
There are many equivalent ways to state that an element is elliptic as can be seen in the following theorem.
\begin{thm}[{\cite[Proposition 3.2, Definition 3.3]{BradyMcCammond}}]
    \label{thm:elliptic_equiv}
    For an element $w$ in an affine Coxeter group $W$, the following are equivalent
    \begin{itemize}
        \item $w$ elliptic, 
        \item $\Mov(w) \subseteq V$ is a linear subspace, and
        \item $\Mov(w)$ contains the zero vector.
    \end{itemize}
\end{thm}

For every vector $\lambda \in V$ there is a Euclidean isometry $t_\lambda$ called a \defn{translation} which sends $x \in E$ to $x + \lambda$.
An element in an affine Coxeter group $W$ is a translation if and only if it is in the kernel $T$ of the projection map $\pi: W \to W_0$ given in Subsection \ref{ssec:affine_coxeter_groups}.
In this case, the set of translations is precisely the set:
\begin{equation}
    T = \set{t_\lambda}{\lambda = \sum_{\alpha \in \Phi} \frac{2\lambda_\alpha}{\gen{\alpha, \alpha}}\alpha, \, \lambda_\alpha \in \Z} .
    \label{eq:translations}
\end{equation}
(see \cite[Definition 1.20]{LewisMcCammondPetersenSchwer}).
This leads to the following proposition.
\begin{prop}[{\cite[Proposition 1.21]{LewisMcCammondPetersenSchwer}}]
    \label{prop:translation_lambda}
    Let $t_\lambda$ be a translation and $w \in W$ an arbitrary element.
    Then $\Mov(t_\lambda w) = \lambda + \Mov(w) =\set{\lambda+m}{m\in\Mov(w)}$.
\end{prop}

In particular, elliptic elements tell us when the product of reflections is linearly independent.
For $r_i \in R$, let $\alpha_i$ be the root associated to $r_i$ and $H_i$ be the hyperplane associated to $\alpha_i$.
\begin{lem}[{\cite[Lemma 1.26]{LewisMcCammondPetersenSchwer}}]
    \label{lem:length_is_ind}
    Let $w = r_1 \ldots r_k$ be a product of reflections in $W$.
    If $w$ is elliptic and $\ell_R(w) = k$ then the $\alpha_i$ are linearly independent.
    Conversely, if the $\alpha_i$ are linearly independent, then $w$ is elliptic, $\ell_R(w) = k$, $\Fix(w) = H_1 \cap \cdots\cap H_k$ and $\Mov(w) = \vspan\left( \left\{ \alpha_1, \ldots, \alpha_k \right\} \right)$.
\end{lem}

\begin{exa}
    \label{ex:affineA2Mov}
    Continuing from Example~\ref{ex:affineA2} let $w = st$ be an element in the affine Coxeter group $\tilde{A}_2$.
    Recall that we can rewrite $st$ as $t(sts)$ and $(sts)s$.
    Since $w$ is in the finite part, the origin is in its fixed space and therefore it is elliptic.
    Furthermore, we recall that $s = r_{\alpha_s, 0}$, $t = r_{\alpha_t, 0}$ and $sts = r_{\alpha_r, 0}$.
    In light of Lemma~\ref{lem:length_is_ind}, we know that
    \[
        \Mov(w) = \vspan\left( \left\{ \alpha_s, \alpha_t \right\} \right) = \vspan\left( \left\{ \alpha_s, \alpha_r \right\} \right) = \vspan\left( \left\{ \alpha_t, \alpha_r \right\} \right).
    \]
    In particular $\Mov(w)$ is equal to the entire vector space and therefore contains the origin.
    The fixed space is precisely the origin and $\ell_R(w) = 2$.

    Let $t_\lambda = rsts$ be a translation as shown in the following figure.
    \begin{center}
        \includegraphics[scale=0.2]{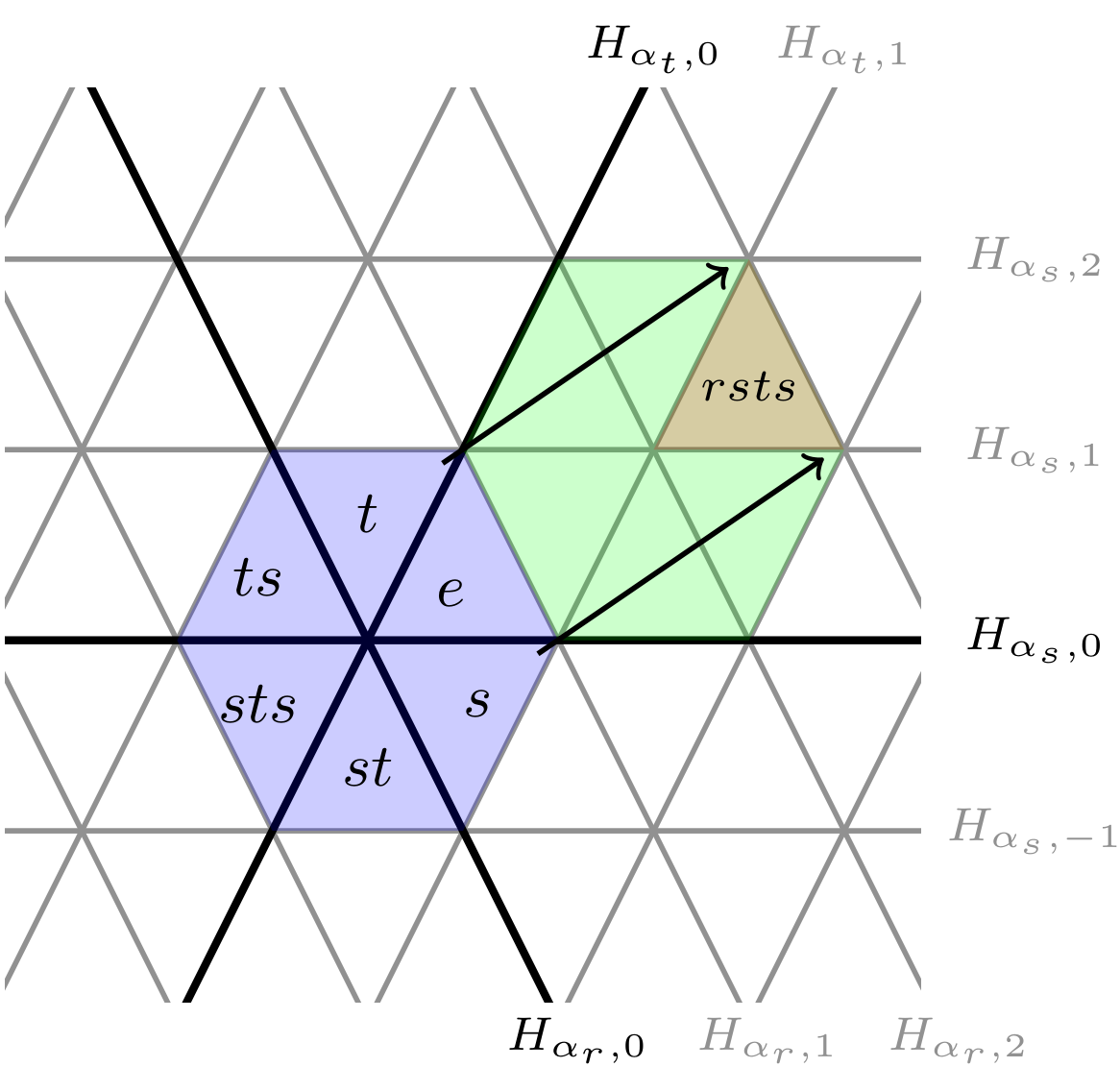}
    \end{center}
    The arrows in the figure represent the vector $\lambda$.
    The translation by $\lambda$ moves the finite part from the blue shaded area in the bottom left to the green and red shaded area in the top right and, in general, moves every point by the vector $\lambda$.

    Consider next $t_\lambda w = rsts\cdot st = rs$.
    Then $t_\lambda w$ is elliptic since
    \[
        \Fix(t_\lambda w) = H_{\alpha_r, 1} \cap H_{\alpha_s, 0} \neq \varnothing.
    \]
    By Lemma~\ref{lem:length_is_ind}, we also have that
    \[
        \Mov(t_\lambda w) = \vspan\left( \left\{ \alpha_s, \alpha_r \right\} \right).
    \]
\end{exa}

\subsection{Dimensions}
\label{ssec:dimensions}
Suppose that $V$ is a Euclidean vector space which contains a (finite) crystallographic root system $\Phi$.
Given a root system $\Phi$, a \defn{root space} is a subspace $U \subseteq V$ such that $U = \vspan\left( U \cap \Phi \right)$.
The \defn{root space arrangement} is the collection of all possible root spaces in $V$.
By construction, as $\Phi$ is finite, there are only a finite number of root spaces in the root space arrangement.
The \defn{root dimension $\dim_{\Phi}(A)$} of a subset $A \subseteq V$ is the minimal dimension of the root spaces contained in the root space arrangement which contains $A$.
This allows us to talk about root dimensions of elements.

The \defn{dimension $\dim(w)$} of an element is the root dimension of its move-set.
In other words, $\dim(w) = \dim_{\Phi}(\Mov(w))$.
Let $\pi$ be the projection map of $W$ to its finite part defined in Subsection~\ref{ssec:affine_coxeter_groups}.
Then let $e(w) = \dim(\pi(w))$ (which we call the \defn{elliptic dimension}) and let $d(w) = \dim(w) - e(w)$ (which we call the \defn{differential dimension}).
The reason for these definitions is that they carry a geometric intuition on what the element $w$ is.

\begin{prop}[{\cite[Proposition~1.31]{LewisMcCammondPetersenSchwer}}]
    \label{prop:dim_and_trans_ellip}
    Let $w$ be an element in an affine Coxeter group.
    \begin{itemize}
        \item Then $w$ is a translation if and only if $e(w) = 0$.
        \item Then $w$ is elliptic if and only if $d(w) = 0$.
    \end{itemize}
\end{prop}

In particular, these dimensions precisely describe the reflection length of an element.
\begin{thm}[{\cite[Theorem A]{LewisMcCammondPetersenSchwer}}]
    \label{thm:A}
    Let $w$ be an element in an affine Coxeter group.
    Then
    \[
        \ell_R(w) = 2 d(w) + e(w)  = 2 \dim(w) - \dim\left( \pi(w) \right)
    \]
    where $d(w)$ and $e(w)$ are defined above.
\end{thm}

The previous theorem implies that if $w$ is in the finite part of $W$, then $\ell_R(w) = e(w)$.
Splitting  $w$ into its differential and elliptic dimensions allows for a factorisation into a translation and an elliptic part $w = t_\lambda u$ which is called a \defn{translation elliptic factorisation of $w$}.
Note that these factorisations are not unique in general.
We will make use of this factorisation when looking at the growth of conjugacy classes in affine Coxeter groups.

Recall from Subsection~\ref{ssec:affine_coxeter_groups} that the projection map $\pi: W \to W_0$ has as a kernel the subgroup of translations $T$,\ie $W_0 \iso W/T$.
As $T$ is finite-index abelian, $W$ is virtually abelian.
In particular, since $T$ is a normal subgroup and $W_0 \cap T = \left\{ 1_W \right\}$ this translation elliptic factorisation implies we have the following semidirect product $W \iso T \rtimes W_0$.
We discuss conjugacy class growth in relation to virtually abelian groups in the next section and use affine Coxeter groups as an example.

\subsection{Conjugacy class growth}
Familiar examples of virtually abelian groups are provided by the affine Coxeter groups which were defined in full detail in Section \ref{sec:coxeter_groups}.
Since every affine Coxeter group $W$ has the form $W \iso T \rtimes W_0$ where $T\iso \Z^d$ is the free abelian normal subgroup of translations and \(W_0\) is the projection of $W$ to its finite part, we have that $T$ is a finite index subgroup and so $W$ is itself virtually abelian.
Since it's also finitely generated, this allows us to use the lemmas from the previous subsection.
In this subsection we give a precise formula for the growth function of a conjugacy class of an arbitrary affine Coxeter group.

\begin{rem}\label{rem:MST2}
    Note that the idea of the proof of the following Lemma and Theorem are similar to that of \cite[Theorem 1.8]{MST}. See also Remark \ref{rem:MST1}.
\end{rem}
\begin{lem}
    \label{lem:rank_is_dim}
    Let \(W = T\rtimes W_0\) be an affine Coxeter group with translations $T$ and finite part $W_0$.
    Let $w \in W_0$ be an arbitrary element in the finite part of $W$.
    Then
    \[
        \ell_R(w) = \dim_{\Phi}\left( \bigcup_{t \in T}\Mov([t,w]) \right) .
    \]
\end{lem}
\begin{proof}
    The set $T$ of translations are linear combinations of roots by (\ref{eq:translations}).
    Note that $ut_{\lambda}u^{-1} = t_{u(\lambda)}$ and $t_\lambda^{-1} = t_{-\lambda}$ for any $u \in W$ and $t_\lambda \in T$.
    Putting these together, we have $\Mov([t_\lambda, w]) = \Mov(t_\lambda t_{w(-\lambda)}) = \left\{\lambda + w(-\lambda)\right\}$ by Proposition~\ref{prop:translation_lambda}.
    Therefore,
    \[
        \dim_{\Phi}\left( \bigcup_{t_\lambda \in T}\Mov([t_\lambda, w]) \right) = \dim_{\Phi}\left(\bigcup_{t_\lambda \in T}\left\{\lambda + w(-\lambda)\right\}\right).
    \]
    Letting $w = r_1 \ldots r_k$ be a reduced expression for $w$ over $R$, then the $\alpha_i \in \Phi$ associated to the $r_i$ are linearly independent by Lemma~\ref{lem:length_is_ind} such that $\Mov(w) = \vspan(\set{\alpha_i}{i \in [k]})$.
    Moreover,
    \begin{align*}
        \Mov(w) &=  \set{\mu \in V}{w(x) = x + \mu \text{ for some } x \in E}\\
        &= \set{\mu \in V}{w(\mu_x) = \mu_x + \mu \text{ for some } \mu_x \in V}\\
        &= \vspan\left( \set{\alpha_i \in \Phi}{w(\lambda) = \lambda + c\alpha_i \text{ for some } t_\lambda \in T,\, i\in [k],\, c \in \R\bs \left\{0\right\}} \right) .
    \end{align*}

    The second equality comes from the fact that any point $x$ is associated with a unique vector $\lambda_x$ where $0 + \lambda_x = x$ (where we have let $w$ denote both the Euclidean isometry on points and the isometry on vectors associated to it).
    The final equality comes from the fact that any isometry in a Euclidean vector space is a linear isomorphism, in addition to the facts that $\vspan\left( \Phi \right) = V$ and $\Mov(w) = \vspan(\set{\alpha_i}{i \in [k]})$.
    Note that for an arbitrary reflection $t_\lambda \in T$ such that $\mu = w(\lambda) - \lambda$ then $w(-\lambda) + \lambda = -w(\lambda) - (-\lambda) = -(w(\lambda) - \lambda) = -\mu$.
    This implies that:
    \begin{align*}
        \dim(w) &=  \dim_{\Phi}\left( \Mov(w) \right)\\
        &= \dim_{\Phi}\left( \vspan\left( \set{\alpha_i \in \Phi}{w(\lambda) = \lambda + c\alpha_i \text{ for some } t_\lambda \in T,\, i\in [k],\, c \in \R\bs \left\{0\right\}} \right) \right)\\
        &= \dim_{\Phi}\left(\bigcup_{t_\lambda \in T}\left\{\lambda + w(-\lambda)\right\}\right).
    \end{align*}
    As $w \in W_0$ is elliptic, then $\dim(w) = d(w) + e(w) = \ell_R(w)$ by Proposition~\ref{prop:dim_and_trans_ellip} and Theorem~\ref{thm:A}.
    Therefore $\ell_R(w) = \dim_{\Phi}\left( \bigcup_{t \in T}\Mov([t,w]) \right)$ as desired.
\end{proof}

\begin{thm}
    \label{thm:finite_part_rank}
    \label{thm:poly_growth_Coxeter}
    Let \(W = T\rtimes W_0\) be an affine Coxeter group with translations $T$ and finite part $W_0$.
    Let $w \in W$ be an arbitrary element in $W$ with translation elliptic factorisation $w = tu$.
    Then the growth rate of the conjugacy class $\cc{w}$ is equivalent to \(n^{\ell_R(u)}\), where \(\ell_R(u)\) is the reflection length.
\end{thm}
\begin{proof}

    Let $w \in W$ be an arbitrary element with translation elliptic factorisation $w = tu$ with $t \in T$ and $u \in W_0$. 
    By Lemma~\ref{lem:vabconj}, since $T$ is a finite-index normal subgroup and since $W_0$ is the set of coset representatives, letting $u \in W_0$ be our fixed representative, we have
    \[
        \cc{w} = \bigcup_{v \in W_0}[T, vuv^{-1}]vwv^{-1}.
    \]
    Then by Lemma~\ref{lem:union} the conjugacy class growth rate of $\cc{w}$ is equivalent to the maximal growth rate of the subsets $[T, vuv^{-1}]vwv^{-1}$ as $v$ ranges across $W_0$.
    Furthermore, by Corollary~\ref{cor:translatedgrowth}, as $vwv^{-1} \in W$, the growth rate of each $[T, vuv^{-1}]vwv^{-1}$ is equivalent to the growth rate of $[T, vuv^{-1}]$.
    In other words, the growth rate of the conjugacy class $\cc{w}$ is equivalent to the maximal growth rate of the subsets $[T, vuv^{-1}]$ as $v$ ranges across $W_0$.

    Since any conjugate of a reflection is again a reflection, the function $\ell_R$ is constant on conjugacy classes and in particular $\ell_R(vuv^{-1}) = \ell_R(u)$.
    Therefore it suffices to show that an arbitrary subset $[T, vuv^{-1}]$ for $v \in W_0$ has growth rate equivalent to $n^{\ell_R(vuv^{-1})} = n^{\ell_R(u)}$.

    Let $u' = vuv^{-1}$ for some $v \in W_0$.
    By Lemma~\ref{lem:commlength}, $[T, u'] = \set{[t,u']}{t \in T}$ as $T$ is a normal abelian subgroup of $W$.
    Furthermore, by normality of $T$, $[T, u']$ is contained in $T$ and is therefore abelian of finite rank.
    By Proposition~\ref{prop:abeliangrowth}, it has polynomial growth whose degree is equal to its rank.
    Therefore it suffices to show $\vrank\left([T, u']\right) = \ell_R(u')$.

    Since the move set of a translation $t \in T$ is precisely the set $\{\lambda\}$ for some $\lambda \in V$, then a translation $t$ is a product of translations $t_1 \ldots t_n$ if and only if $\lambda = \lambda_1 + \cdots + \lambda_n$.
    Therefore, the rank of $[T, u'] \subseteq T$ is equal to the dimension of the span of the move sets of the elements of $[T, u']$\:
    \[
        \vrank([T, u']) = \dim(\vspan(\bigcup_{t \in T} \Mov([t, u']))).
    \]
    Since $[T, u'] \subseteq T$, the elements of  $[T, u']$ are linear combinations of roots by (\ref{eq:translations}).
    This implies that $\vspan(\bigcup_{t \in T} \Mov([t, u']))$ is a root space.
    In other words
    \[
        \vrank\left( [T, u'] \right) = \dim_{\Phi}\left( \bigcup_{t \in T} \Mov\left( [t, u'] \right) \right).
    \]
    By Lemma~\ref{lem:rank_is_dim}, since $u' \in W_0$, then $\dim_{\Phi}\left( \bigcup_{t \in T} \Mov\left( [t, u'] \right) \right) = \ell_R(u')$, as desired.
\end{proof}

\section*{Acknowledgements}
The second named author thanks Denis Osin for the initial suggestion of studying the growth of a conjugacy class.
The authors would also like to thank the reviewer(s) for helpful suggestions on improving this article.


\bibliography{references}{}
\bibliographystyle{plain}

\bigskip

\end{document}